\documentclass[11pt]{amsart}

\usepackage{amssymb, amscd}
\usepackage{graphicx}

\newcommand\Dt{\Delta}

\newcommand\ld{\lambda}

\newcommand\x{\xi}

\newcommand\Ph{\Phi}

\newcommand\ph{\varphi}

\newcommand\CC{\mathbb C}

\newcommand\ZZ{\mathbb Z}

\newcommand\SU{{\mathrm{SU}}}
\newcommand\U{{\mathrm{U}}}

\newcommand\End{\operatorname{End}}

\makeatletter 
\newcommand\matc[4]{\left( {#1\@@atop #3}{#2\@@atop #4}\right)}

\newcommand\matr[4]{\left( {\hfill #1\@@atop\hfill #3}{\hfill
#2\@@atop\hfill #4}\right)}

\newcommand\matl[4]{\left( { #1\@@atop #3}{ #2\@@atop\hfill #4}\right)}

\makeatother

\newcommand\widearray[1]{\renewcommand\arraystretch{1.4} \begin{array}{#1}}

\newcommand\lw[1]{\hbox{}_{#1}\!}

\newcommand\vzm[2]{\mathchoice{\left\{\, #1 : #2 \,\right\}}{\{\, #1
:\allowbreak #2 \,\}}{\{ #1 :\allowbreak #2 \}} {\{ #1 :\allowbreak #2 \}}}

 \theoremstyle{plain}
\newtheorem{thm}{Theorem}[section]

\newtheorem*{prop4.5}{Proposition \ref{expD2}}
\newtheorem*{prop4.6}{Proposition \ref{expE1}}
\newtheorem*{prop4.7}{Proposition \ref{expE2}}
\newtheorem*{lem5.3}{Lemma \ref{dderx1}}
\newtheorem*{lem5.4}{Lemma \ref{ddery1}}
\newtheorem*{lem5.5}{Lemma \ref{ddery2}}
\newtheorem*{lem5.6}{Lemma \ref{dderx2}}
\newtheorem*{lem5.7}{Lemma \ref{ddery1x1}}
\newtheorem*{lem5.8}{Lemma \ref{derx2y1}}
\newtheorem*{lem5.9}{Lemma \ref{derx2y2}}

\newtheorem{conj}[thm]{Conjecture}

\theoremstyle{definition}
\newtheorem{defn}[thm]{Definition}

\theoremstyle{remark}

\title[matrix valued spherical functions]{\sc An invitation to matrix
valued spherical
functions:
linearization of products in the case of the complex projective space
$P_2(\CC)$}

\author{ F. A. Gr\"unbaum}
\address{Departament of Mathematics, University of California, Ber\-ke\-ley
CA 94705}
\email{grunbaum@math.berkeley.edu}

\author{I. Pacharoni}
\address{CIEM-FaMAF, Universidad Nacional de C\'or\-do\-ba,
C\'or\-do\-ba~5000, Argentina}
\email{pacharon@mate.uncor.edu}

\author{J. Tirao}
\address{CIEM-FaMAF, Universidad Nacional de C\'or\-do\-ba,
C\'or\-do\-ba~5000, Argentina}
\email{tirao@mate.uncor.edu}

\thanks{This paper is partially supported by NSF grants FD9971151 and
1-443964-21160
and by CONICET grant PIP 655-98}
\begin{document}
\maketitle

\begin{abstract}

The classical (scalar valued) theory of spherical functions (put forward
by Cartan and others after him) allows one to unify under one roof a number
of
examples that were very well known before the theory was formulated.
These examples include many special functions like Jacobi polynomials,
Bessel functions, Laguerre polynomials, Hermite polynomials,
Legendre functions, etc.

All these functions had ``proved themselves'' as the work-horse in many
areas of mathematical physics before the appearance of a unifying theory.
Many of these functions have found interesting applications in signal
processing in general as well as in very specific areas like medical
imaging.

The theory of matrix valued spherical functions
gives a natural
extension of the well-known theory for the scalar valued case.

The historical development in the matrix valued extension of this theory is
entirely different. In this case the theory has gone ahead of the examples.

The purpose of this note is to give some pointers to some examples
and to "tease" some readers into this new aspect in the world of special
functions.

We close with a remark connecting the functions described here with the
theory of matrix valued orthogonal polynomials.

\end{abstract}
\maketitle

\section{Introduction and statement of results}

The theory of matrix valued spherical functions, see \cite{GV} and \cite{T},
gives a natural
extension of the well-known theory for the scalar valued case, see
\cite{He}. We start with a few remarks about the scalar valued case.

The classical (scalar valued) theory of spherical functions (put forward
by Cartan and others after him) allows one to unify under one roof a number
of
examples that were very well known before the theory was formulated.
These examples include many special functions like Jacobi polynomials,
Bessel functions, Laguerre polynomials, Hermite polynomials,
Legendre functions, etc.

All these functions had ``proved themselves'' as the work-horse in many
areas of mathematical physics before the appearance of a unifying theory.
Many of these functions have found interesting applications in signal
processing in general as well as in very specific areas like medical
imaging.
It suffices to recall, for instance, that Cormack's approach
\cite{C}---for which
he got the 1979 Nobel Prize in Medicine, along with G.\
Hounsfield---was based on
classical orthogonal polynomials and that the  work of Hammaker
and Solmon \cite{HS}  as well as that of Logan and Shepp \cite{LS}
is based on the use of Chebychev polynomials.

The crucial property
here is the fact that these functions satisfy the integral equation that
characterizes  spherical functions of a
homogeneous space.
For a review on some of these topics the reader can either look at 
some of the specialized books on the subject such as \cite{He}
or start from a more
introductory approach as that given in either \cite{DMcK} and \cite{T1}.

This integral equation is actually satisfied by all
Gegenbauer polynomials and not only those corresponding to
symmetric spaces. This point is fully exploited in
  \cite{DG} where this property is put to use to show that different weight
functions can be used in carrying out the usual tomographic
operations of projection and
backprojection. This works well for parallel beam tomography but has never
been made to work for fan beam tomography because of a lack of an underlying
group theoretical formulation in this case. For a number of issues in this
area, including a number of open problems, see \cite{G2}.

For a variety of other applications of spherical functions
one can look at \cite{DMcK} , \cite{T1},
\cite{T2}.

We now come to the main issue in this note.

The situation with the matrix valued extension of this theory is
entirely different. In this case the theory has gone ahead of the examples
and, in fact, to the best of our knowledge, the first examples
involving non-scalar matrices
have been given recently in \cite{GPT1} , \cite{GPT2} and \cite{GPT3}.
For scalar
valued instances
of non-trivial type, see \cite{HeSc}.

The issue of how useful these functions may turn out to be as a tool
in areas like geometry, mathematical physics, or signal processing in the
broad sense is still open. From a historical perspective
one could argue, rather tautologically, that the usefulness of the classical
spherical functions rests on the many interesting properties they all share.
With that goal in mind, it is natural to try to give a glimpse at these new
objects and to illustrate some of their properties. The rather mixed
character
of the audience attending these lectures gives us an extra incentive to make
this material accessible to people that might normally not look in the
specialized literature.

The purpose of this contribution is thus
to present  very briefly the essentials of the theory and
to describe one example in some detail. This is not the appropriate place
for a complete description, and we refer the interested reader to the
papers \cite{GPT1} , \cite{GPT2} and \cite{GPT3}.

We hope to pique the curiosity of some readers by exploring
the extent to which the property
of ``positive linearization of products'' holds
in the case of the spherical functions
associated to 
$P_2(\CC)$.
This result has been important in the scalar case,
including its use in the proof of the Bieberbach conjecture, see \cite{AAR}.
The property in question is illustrated well by considering the case of
Legendre polynomials: the product of any two such is expressed as a linear
combination involving other Legendre polynomials with degrees ranging from
the absolute value of the difference to the sum of the degrees of the two
factors involved. Moreover, the coefficients in this expansion are positive.

We should stress that the intriguing property described here
is one enjoyed by a {\em matrix valued function} put together from different
spherical
functions of a given type. In the classical {\em scalar valued} case
these two notions
agree and the warning is not needed.
This combination of spherical functions has already been seen, see
\cite{GPT1} , \cite{GPT2} and \cite{GPT3}  to enjoy a natural form
of the {\em bispectral
property}. For an introduction to this expanding subject we could 
consult, for instance, \cite{DG1}, \cite{G12}. The roots of this problem
are too long to trace in this short paper, but the reader may want
to take a look at \cite{S1}. For off-shoots that have yet to be explored
further one can also see \cite{G13} and \cite{G15}. The short version of
the story is that some remarkably useful {\em algebraic} properties that have
surfaced first in signal processing and which one would like to extend
and better understand have a long series of connections with other parts
of mathematics. For a collection of problems arising in this area see
\cite {HK}.

The issue of linearization of products, without insisting on any positivity
results, plays (in the scalar valued case)  an important role
in fairly successful
applications of mathematics. For example,
the issue of expressing the product of spherical
harmonics of different degrees as a sum of spherical
harmonics plays a substantial role in both theoretical
and practical algorithms for the harmonic analysis of functions on the
sphere.
For some developments in this area see \cite{DH} as well as \cite{KMHR}.

In the context of Quantum Mechanics this discussion is the backbone of the
{\em addition rule for angular momenta} as can be seen in any textbook on
the subject.

In the last section we make a brief remark connecting the functions described
here with the theory of matrix valued orthogonal polynomials, as developed for
instance in \cite{D} and  \cite{DVA}.

\section{Matrix valued spherical functions}

Let $G$ be a locally compact unimodular group and let $K$ be a compact
subgroup of $G$. Let $\hat K$ denote the set of all
equivalence classes of complex finite dimensional irreducible
representations of $K$; for each $\delta\in \hat K$, let
$\x_\delta$ denote the character of $\delta$, $d(\delta)$ the degree of
$\delta$, i.e.\ the dimension of any representation in
the class $\delta$, and $\chi_\delta=d(\delta)\x_\delta$.

Given a homogeneous space $G/K$ a zonal spherical function
(\cite{He}) $\ph$ on $G$ is a
continuous complex valued function which satisfies
$\ph(e)=1$ and
\begin{equation}\label{defclasica}
\ph(x)\ph(y)=\int_K \ph(xky)\, dk \qquad \qquad x,y\in G.
\end{equation}

A fruitful generalization of the above concept is given in the following
definition.

\begin{defn}[\cite{T},\cite{GV}]\label{defesf}
A spherical function $\Ph$ on $G$ of type $\delta\in \hat K$ is a
continuous function on $G$ with values in $\End(V)$ such
that
\begin{enumerate}
\item[i)] $\Ph(e)=I$ ($I$= identity transformation).

\item[ii)] $\Ph(x)\Ph(y)=\int_K \chi_{\delta}(k^{-1})\Ph(xky)\, dk$, for
all $x,y\in G$.
\end{enumerate}
\end{defn}

The connection with differential equations of the group $G$ comes about from
the property below.

Let $D(G)^K$ denote the algebra of all
left invariant differential operators on $G$ which are also invariant under
all
right translation by elements in $K$.
If $(V,\pi)$ is a finite dimensional
irreducible representation of $K$ in the equivalence
class $\delta \in \hat K$, a spherical function on $G$ of
type $\delta$ is characterized by
\begin{enumerate}
\item[i)] $\Ph:G\longrightarrow \End(V)$ is analytic.
\item[ii)] $\Ph(k_1gk_2)=\pi(k_1)\Ph(g)\pi(k_2)$, for all $k_1,k_2\in K$,
$g\in G$, and $\Phi(e)=I$.
\item[iii)] $[D\Ph ](g)=\Ph(g)[D\Ph](e)$, for all $D\in D(G)^K$, $g\in G$.
\end{enumerate}
\smallskip
We will be interested in the specific example
given by the complex projective plane. This can be realized
as the homogeneous space $G/K$, where $G=\SU(3)$ and
$K={\mathrm{S}}(\U(2)\times\U(1))$. In this case iii) above can be replaced
by:
$[\Dt_2\Ph ](g)=\ld_2\Ph(g)$, $[\Dt_3\Ph ](g)=\ld_3\Ph(g)$ for
all $g\in G$ and for some $\ld_2,\ld_3\in \CC$.
Here $\Dt_2$ and $\Dt_3$ are two algebraically independent generators of the
polynomial algebra $D(G)^G$ of all differential operators on $G$ which are
invariant under left and right multiplication by elements in $G$.

The set $\hat K$ can be identified with the set $\ZZ\times \ZZ_{\geq 0}$.
If $k=\left(\begin{smallmatrix} A&0\\ 0& a
\end{smallmatrix}\right)$, with $A\in \U(2)$ and $a=(\det A)^{-1}$, then
$$\pi(k)=\pi_{n, l}(A)=\left(\det A\right )^n\,A^ l,$$
where $A^ l$ denotes the $ l$-symmetric power of $A$, defines an
irreducible representation of $K$ in the class
$(n, l)\in \ZZ\times \ZZ_{\geq 0}$.

For simplicity we restrict ourselves in this brief presentation
to the case $n\geq 0$. The paper \cite{GPT1} deals with the general case.
The
representation $\pi_{n, l}$ of $\U(2)$ extends to a unique holomorphic
multiplicative map of
${\mathrm M}(2,\CC)$ into $\End(V_\pi)$, which we shall still denote by
$\pi_{n, l}$. For any
$g\in {\mathrm M}(3,\CC)$, we shall denote by $A(g)$ the left upper
$2\times 2$ block of $g$, i.e.
$$A(g)=\matc{g_{11}}{g_{12}}{g_{21}}{g_{22}}.$$

For any $\pi=\pi_{(n, l)}$ with $n\geq 0$ let $\Phi_\pi:G\longrightarrow
\End(V_\pi)$ be defined by
$$\Phi_\pi(g)=\Phi_{n,l}(g)= \pi_{n, l}(A(g)).$$
It happens that $\Phi_\pi$ is a spherical
function of type $(n, l)$, one that will play a very important role in
the construction of {\em all} the remaining spherical functions of the same
type.

Consider the
open set
$${\mathcal A}=\vzm{g\in G}{\det A(g)\neq 0}.$$
The group $G=\SU(3)$ acts in a natural way in the complex projective plane
$P_2(\CC)$. This action is transitive and $K$ is
the isotropy subgroup of the point $(0,0,1)\in P_2(\CC)$. Therefore $
P_2(\CC)= G/ K$. We shall identify the complex plane
$\CC^2$ with the affine plane $\vzm{(x,y,1)\in P_2(\CC)}{(x,y)\in \CC^2}$.

The canonical projection $p:G\longrightarrow P_2(\CC)$ maps the open dense
subset ${\mathcal A}$ onto the affine plane
$\CC^2$. Observe that ${\mathcal A}$ is stable by left and right
multiplication by elements in $K$.

To determine all spherical functions
$\Phi:G\longrightarrow \End(V_\pi)$ of type $\pi=\pi_{n,l}$, we use the
function $\Phi_\pi$ introduced above
in the following way: in the open
set ${\mathcal A}$ we define a function $H$ by
$$H(g)=\Phi(g)\, \Phi_\pi(g)^{-1},$$
where $\Phi$ is suppose to be a spherical function of type $\pi$.
Then $H$ satisfies
\begin{enumerate}
\item [i)] $H(e)=I$.
\item [ii)]$ H(gk)=H(g)$, for all $g\in {\mathcal A}, k\in K$. \item [iii)]
$H(kg)=\pi(k)H(g)\pi(k^{-1})$, for all
$g\in {\mathcal A}, k\in K$.
\end{enumerate}

\noindent Property ii) says that $H$ may be considered as a function on
$\CC^2$.

The fact that $\Phi$ is an eigenfunction of $\Delta_2$ and $\Delta_3$
makes $H$ into an eigenfunction of certain
differential operators $D$ and $E$ on $\CC^2$.

We are interested in considering the differential operators $D$ and $E$
applied to a function $H\in C^\infty(\CC^2)\otimes
\End(V_\pi)$ such that $H(kp)=\pi(k)H(p)\pi(k)^{-1}$, for all $k\in K$ and
$p$ in the affine complex plane $\CC^2$. This
property of $H$ allows us to find ordinary differential operators $\tilde
D$ and $\tilde E$ defined on the interval
$(0,\infty)$ such that $$ (D\,H)(r,0)=(\tilde D\tilde H)(r),\qquad
(E\,H)(r,0)=(\tilde E\tilde H)(r),$$ where $\tilde
H(r)=H(r,0)$.

Introduce the variable $t=(1+r^2)^{-1}$ which converts the
operators $\tilde D$ and $\tilde E$
into new operators $D$ and $E$.

The functions $\tilde H$ turn out to
be diagonalizable. Thus, in an appropriate
basis of
$V_\pi$, we can write $\tilde H(r) = H(t) = (h_0(t),\cdots,
h_ l(t))$.

We find it very convenient to introduce
two integer parameters $w$, $k$ subject to the following three inequalities:
$0\leq w$, $0\leq k\leq  l$,
which give a very convenient parametrization of the irreducible spherical
functions of type $(n, l)$. In fact, for each pair $( l,n)$, there are
a total of $l+1$ {\em families} of matrix valued functions of $t$ and $w$.
In this instance these matrices are diagonal and one can put these diagonals
together into a {\em full matrix valued function} as we will do in the next
two sections.
It appears that this function, which concides with the usual spherical
function
in the scalar case, enjoys some interesting properties.

 The reader can consult \cite{GPT1} to find a fairly detailed description
of the entries that make up the matrices mentioned up to now. A flavor
of the results is given by the following statement.

For a given $l \ge 0$, the spherical functions corresponding to the pair
$(l,n)$ have components that are expressed in
terms of generalized hypergeometric functions of the form
$\lw{p+2}F_{p+1}$, namely
\[
\lw{p+2}F_{p+1} \left( \begin{smallmatrix} a,b,s_1+1,\dots,s_p+1
\vspace{2pt} \\
c,s_1,s_2,\dots,s_p
\end{smallmatrix} ;t \right) = \sum_{j=0}^{\infty} \frac
{(a)_j(b)_j}{j!(c)_j} (1 + d_1j + \dots + d_pj^p)t^j.
\]

\section{The bispectral property}

For given non-negative integers $n$, $l$ and $w$ consider the matrix
whose rows are given by the vectors $H(t)$
corresponding to the values $k=0,1,2,\dots,l$ discussed above. Denote
the corresponding matrix
by
\[ \Phi (w,t). \]

As a function of $t$, $\Phi(w,t)$ satisfies two differential equations
\[
D \Phi(w,t)^t = \Phi(w,t)^t\Lambda , \quad E \Phi(w,t)^t=\Phi(w,t)^tM \ .
\]
Here $\Lambda$ and $M$ are diagonal matrices with $\Lambda (i,i)=
-w(w+n+i+ l +1)- (i-1)(n+i)$,
$M(i,i)=\Lambda(i,i)(n- l+3i-3)-3(i-1)(l-i+2)(n+i)$, $1\leq
i\leq l+1$;
$D$ and $E$ are the differential operators introduced earlier. Moreover we
have

\begin{thm} There exist matrices $A_w$, $B_w$, $C_w$, independent of $t$,
such that
\[
A_w\Phi(w-1,t)+B_w\Phi(w,t)+C_w \Phi(w+1,t) = t\Phi(w,t)\ .
\]
\end{thm}

The matrices $A_w$ and $C_w$ consist of two diagonals each and $B_w$ is
tridiagonal.
Assume, for convenience, that these vectors are normalized in such a way
that for $t=1$ the matrix
$\Phi(w,1)$ consists of all ones.

For details on these matrices as well as for a full proof of this statement,
which was conjectured in \cite{GPT1},
the reader can consult \cite{GPT2}.

\section{Linearization of products}

The property in question states that the product of members of certain
families of (scalar valued) orthogonal polynomials is given by
an expansion of the form
\[
P_{i} P_{j} = \sum_{k=|j-i|}^{j+i} a_{k} P_{k}
\]
and that
the coefficients in the expansion
are all non-negative.

For a nice and detailed
account of the situation in the scalar case, see for instance \cite{A},
\cite{S}. Very important contributions on these and related matters
are \cite{G} and \cite{K}.

It is important to note that the property in question is not true for all
families of orthogonal polynomials, in fact it is not even true for all
Jacobi polynomials $P_w^{(\alpha,\beta)}$, {\em normalized
by the condition $P_w^{(\alpha,\beta)}(1)$ positive}.
For our purpose it is important to recall
that non-negativity is satisfied if  $\alpha \ge \beta$ and
$\alpha+\beta \ge 1$.

\smallskip
\noindent {\bf The case $l=0$, $n>1$.}

 From \cite{GPT1} we know that when $l=0$ and $n \ge 0$ the appropriate
eigenfunctions (without the standard normalization)  are given by
\[
\Phi(w,t) = \lw{2}F_1 \left( \begin{smallmatrix} -w,\,\;w+n+2
\vspace{2pt} \\ n+1
\end{smallmatrix} ;t\right).
\]
This means that with the usual convention that the Jacobi polynomials
are positive for $t=1$ we are dealing with the family
\[
P_w^{(1,n)} (t).
\]

If $n=0$ or $n=1$ the family $P_w^{(1,n)}$ meets
the sufficient conditions for non-negativity given above.
For $n=0$ the coefficients $a_k$ are all strictly positive; in the case
$n=1$ the 
coefficients $a_{|i-j|+{2k}}$, are strictly positive
while the coefficients
$a_{|i-j|+{k}}$, {k} odd,  
are zero, as the example below illustrates.

We now turn our attention to the case $n>1$.

\begin{conj}
For $n$ an integer larger than one,
the coefficients in the expansion for the product $P_{i} P_{j}$  above
alternate in sign.
\end{conj}

This conjecture is backed up by extensive experiments, one of which is shown
below. It deals with the case of w (i.e. i and j ) equal to 3 and 4.
Prof. Richard Askey, supplied a proof of this conjecture. This gives
us a new chance to thank him for many years of encouragement and help.

\medskip
The product of the (scalar valued, and properly normalized) functions
$\Phi(3,t)$  and $\Phi(4,t)$
is given by the expansion
\begin{align*}
\Phi(3,t)\Phi(4,t)  =& \;a_1 \Phi(1,t) + a_2 \Phi(2,t)  + a_3\Phi(3,t) + a_4
\Phi(4,t) \\
& + a_5\Phi(5,t) + a_6   \Phi(6,t)  + a_7  \Phi(7,t)
\end{align*}
with coefficients given by the expressions

\begin{eqnarray*}
a_1 &= &\textstyle{\frac {(n+2)(n+3)(n+4)}{(n+8)(n+9)(n+10)}}
\vspace{2\jot} \\
a_2 &= &\textstyle{-\frac
{6(n-1)(n+3)(n+4)(n+6)^2}{(n+7)(n+8)(n+9)(n+10)(n+11)}}
\vspace{2\jot} \\
a_3 &= &\textstyle{\frac
{3(n+4)(n+5)(7n^3+52n^2+67n+162)}{(n+7)(n+9)(n+10)(n+11)(n+12)}}
\vspace{2\jot} \\
a_4 &= &\textstyle{-\frac
{4(n-1)(n+6)(11n^3+123n^2+436n+648)}{(n+8)(n+9)(n+11)(n+12)(
n+13)}}
\vspace{2\jot} \\
a_5 &= &\textstyle{\frac
{3(n+5)(n+6)(n+7)(19n^3+155n^2+162n+504)}{(n+8)(n+9)(n+10)
(n+11)(n+13)(n+14)}}
\vspace{2\jot} \\
a_6 &= &\textstyle{-\frac
{42(n-1)(n+5)(n+6)^2(n+7)(n+8)}{(n+9)(n+10)(n+11)(n+12)(n+
13)(n+15)}}
\vspace{2\jot} \\
a_7 &= &\textstyle{\frac
{14(n+5)(n+6)^2(n+7)^2(n+8)}{(n+10)(n+11)(n+12)(n+13)(n+14)(n
+15)}}
\end{eqnarray*}

This shows that even in the scalar valued case, as soon as we are dealing
with non-classical spherical functions we encounter an interesting sign
alternating property that is quite different from the more familiar case.
Here and  below we see that things become different once $n$ is an integer
larger than one.

\medskip
Now we explore the picture in the case of general $l$ .

\smallskip
\noindent {\bf The case $l>0$, $n>1$.}

\begin{conj}
If $i \le j$ then the product of $\Phi(i,t)$ and
$\Phi(j,t)$ allows for a (unique)
expansion of the form
\[
\Phi(i,t) \Phi(j,t) = \sum_{k=\min\{j-i-l,0\}}^{j+i+l} A_{k} \Phi(k,t).
\]
Here the coefficients $A_{k}$ are matrices and the matrix valued
function $\Phi(w,t)$ is the one introduced in section $3$. This conjecture
holds for all nonnegative $n$ and is well known for $l=0$ and $n=0$.
\end{conj}

One should remark that in the case of $l=0$ we obtain the usual range in the
expansion coefficients ranging from $j-i$ to $j+i$ as in the case of
addition of angular momenta.
For larger values of $l$ we see that {\em extra}
terms appear at each end of the expansion.

\begin{conj}
If $i<j$ then the coefficients $A_{k}$ in the expansion
\[
\Phi(i,t) \Phi(j,t) = \sum_{k=\min\{j-i-l,0\}}^{j+i+l}  A_{k} \Phi(k,t).
\]
with $k$ in the range $j-i,j+i$ have what we propose to call
``the hook alternating property.''
\end{conj}

We will explain this conjecture by displaying one example. First notice
that we exclude those coefficients that are not in the {\em traditional}
or {\em usual} range discussed above.

At this point it may be appropriate in the name of {\em truth in advertisement}
to admit that we have no concrete evidence of the significance of the property
alluded to above and displayed towards the end of the paper. We trust
that the reader will find the property cute and intriguing. 
It would be
very dissapointing if nobody were to find some use for it.  

The results illustrated below have been checked for many values of $l>0$,
but are displayed here for $l=1$ only.

\medskip
Recall that from \cite{GPT1} the rows that make up the matrix valued
function $H(t,w)$ are given as follows:
the first row is obtained from the column vector
\[
H(t) = \left( \begin{array}{c}
\left(1-\frac{\lambda}{n+1} \right) \lw{3}F_2 \left( \begin{smallmatrix}
-w,\,\;w+n+3,\,\;\lambda-n \vspace{2pt} \\
n+2,\,\;\lambda - n-1
\end{smallmatrix} ;t \right) \vspace{2pt} \\ \lw{2}F_1
\left(\begin{smallmatrix}
-w,\,\;w+n+3 \vspace{2pt} \\
n+1
\end{smallmatrix} ;t \right)
\end{array} \right)
\]
with \[\lambda =-w(w+n+3)\]
and the second row comes from the column vector
\[
H(t) = \left( \begin{array}{c}
\lw{2}F_1 \left( \begin{smallmatrix}
-w,\,\;w+n+4 \vspace{2pt} \\
n+2
\end{smallmatrix} ;t\right) \vspace{2pt} \\ -(n+1)\, \lw{3}F_2
\left(\begin{smallmatrix} -w-1,\,\; w+n+3,\,\;\lambda
\vspace{2pt} \\ n+1,\,\;\lambda-1
\end{smallmatrix} ;t\right)
\end{array} \right)
\]
with
\[
\lambda = -w(w+n+4)-n-2.
\]

\noindent The product of the matrices $\Phi(2,t)$ and $\Phi(6,t)$  is given
by
the expansion
\begin{align*}
\Phi(2,t) \Phi(6,t) = &\; A_3 \Phi(3,t)  + A_4 \Phi(4,t)  + A_5
\Phi(5,t) + A_6 \Phi(6,t)  \\
 &\;+  A_7 \Phi(7,t) +
A_8 \Phi(8,t)  + A_9 \Phi(9,t)
\end{align*}
where the coefficient matrices are displayed below.
\noindent The matrix
\[
A_3=\begin{pmatrix}
0 & 0 \\
0 & \frac {16(n+4)(n+5)(n+6)^2(n+7)^2}{(n+11)(n+12)(n+13)(n+14)(n+15)(n+16)}
\end{pmatrix}.
\]

\medskip
\noindent The matrix
\[
A_4=\matc{L_{11}}{L_{12}}{L_{21}}{L_{22}}
\]
with
\begin{align*}
L_{11} &=\tfrac {15(n+5)^2(n+6)(n+8)}{2(n+12)(n+13)(n+14)(n+15)},
\displaybreak[0]\\
L_{12}&= \tfrac
{5(n+5)(n+6)(4\,n^2+55\,n+216)}{6(n+13)(n+14)(n+15)(n+16)},\displaybreak[0]\
\
L_{21}&= \tfrac
{(n+5)(n+6)(n+7)(8\,n^2+153\,n+724)}{2(n+12)(n+13)(n+14)(n+15)(n+16)},
\displaybreak[0]\\
L_{22}&=-\tfrac
{5(n+6)(n+7)(248\,n^4+4665\,n^3+27202\,n^2+45137\,n-23252)}{12(n+11)(n+13)(n
+14)(n+15)(n+16)(n+17)}.
\end{align*}

\medskip
\noindent The matrix
\[
A_5=\matc{M_{11}}{M_{12}}{M_{21}}{M_{22}}
\]
with
\begin{align*}
M_{11} &= \textstyle{-\frac {(n+5)(n+6)(185\,n^3 + 3284\,n^2 + 15732\,n +
10368)}{6(n+7)(n+12)(n+14)(n+15)(n+16)}},\displaybreak[0] \\
M_{12} &= \textstyle{-\frac {(n+5)(85\,n^4 + 1817\,n^3 + 11380\,n^2 +
7072\,n -93460)}{7(n+7)(n+13)(n+15)(n+16)(n+17)}}, \displaybreak[0] \\
M_{21} &=
\textstyle{-\frac {(n+6)^2(170\,n^4 + 4735\,n^3 + 42068\,n^2 + 99767\,n -
168628)}{12(n+7)(n+12)(n+14)(n+15)(n+16)(n+17)}},
\end{align*}
and finally
\begin{align*}
M_{22} = &
\,\textstyle{\frac {4327\,n^7 + 163698\,n^6 + 2480127\,n^5 + 19091004\,n^4
+ 78090428\,n^3
+ 163454544\,n^2}{14(n+7)(n+12)(n+13)(n+15)(n+16)(n+17)(n+18)}}\\
&+\textstyle{\frac { 172290528\,n +
132098688}{14(n+7)(n+12)(n+13)(n+15)(n+16)(n+17)(n+18)}}.
\end{align*}

\medskip
\noindent The matrix
\[
A_6=\matc{N_{11}}{N_{12}}{N_{21}}{N_{22}}
\]
with
\begin{align*}
N_{11} &= \textstyle{\frac {2(193\,n^5 + 5832\,n^4 + 65284\,n^3 +
328884\,n^2 + 727621\,n +
634422)}{7(n+8)(n+13)(n+14)(n+16)(n+17)}} \vspace{2\jot}, \displaybreak[0]\\
N_{12} &= \textstyle{\frac {171\,n^5 + 4729\,n^4 + 45764\,n^3 + 188570\,n^2
+ 442336\,n + 1133640}{8(n+8)(n+14)(n+15)(n+17)(n+18)}} \vspace{2\jot},
\displaybreak[0]\\
N_{21} &= \textstyle{\frac {171\,n^6 + 7071\,n^5 + 116213\,n^4 +
959879\,n^3 + 4245034\,n^2 + 10640548\,n +
15755112}{7(n+8)(n+13)(n+14)(n+16)(n+17)(n+18)}},
\end{align*}
and finally
\begin{align*}
N_{22} =\,&
\textstyle{-\frac {4269\,n^7 + 169934\,n^6 + 2677678\,n^5 + 21066480\,n^4 +
85737209\,n^3 +
169428298\,n^2 }{8(n+8)(n+13)(n+14)(n+15)(n+17)(n+18)(n+19)}}\\
\,&\textstyle{-\frac { 129986220\,n -46794888}
{8(n+8)(n+13)(n+14)(n+15)(n+17)(n+18)(n+19)}}.
\end{align*}

\medskip
\noindent The matrix
\[
A_7=\matc{P_{11}}{P_{12}}{P_{21}}{P_{22}}
\]
with
\begin{align*}
P_{11} &= \textstyle{-\frac {3(n+5)(129\,n^4 + 3710\,n^3  +
36430\,n^2 + 129960\,n + 76536)}{8(n+9)(n+14)(n+15)(n+16)(n+18)}}
\vspace{2\jot}, \displaybreak[0]\\
P_{12} &= \textstyle{-\frac {(n+5)(n+10)(57\,n^3 + 917\,n^2 + 2274\,n -
11268)}{3(n+9)(n+15)(n+16)(n+17)(n+19)}} \vspace{2\jot}, \displaybreak[0]\\
P_{21} &= \textstyle{\frac {-3(57\,n^6 + 2505\,n^5 + 44489\,n^4 +
389955\,n^3 + 1576582\,n^2 + 1465908\,n -4434696)}
{8(n+9)(n+14)(n+15)(n+16)(n+18)(n+19)}},
\end{align*}
and finally
\[
P_{22} =
\textstyle{\frac {2(n+10)(829\,n^6 + 27979\,n^5 + 352571\,n^4 + 2024521\,n^3
+ 5197384\,n^2 + 5712396\,n + 5004720)}
{3(n+9)(n+14)(n+15)(n+16)(n+17)(n+19)(n+20)}}.
\]

\medskip
\noindent The matrix
\[
A_8=\matc{Q_{11}}{Q_{12}}{Q_{21}}{Q_{22}}
\]
with
\begin{align*}
Q_{11} &= \tfrac {5(n+5)(n+6)(21\,n^2 + 401\,n +
1920)}{6(n+15)(n+16)(n+17)(n+18)}, \displaybreak[0]\\
Q_{12} &=  \tfrac {15(n+5)(n+6)(n+8)(n+11)}{2(n+16)(n+17)(n+18)(n+19)},\\
Q_{21} &=   \tfrac {5(n+6)(10\,n^4 + 329\,n^3 + 4942\,n^2 + 36611\,n +
96300)}{6(n+15)(n+16)(n+17)(n+18)(n+20)},\displaybreak[0] \\
Q_{22} &= -\tfrac {3(n+6)(n+11)(430\,n^4 + 9773\,n^3 + 67728\,n^2 +
129129\,n -
59220)}{4(n+15)(n+16)(n+17)(n+18)(n+19)(n+21)}.
\end{align*}

\medskip
\noindent The matrix $A_9$
\[
\begin{pmatrix}
0 & 0 \vspace{2\jot} \\
\frac {99(n+4)(n+6)(n+7)(n+10)}{4(n+16)(n+17)(n+18)(n+19)(n+20)} &
\frac
{165(n+4)(n+6)(n+7)(n+8)(n+10)(n+12)}{2(n+16)(n+17)(n+18)(n+19)(n+20)(n+21)}
\end{pmatrix}.
\]


\medskip
Notice that if we concentrate our attention on the coefficients within the
{\em traditional} range we see that the first matrix $A_4$ has
its first {\em hook}
made up of positive entries, the second hook (which in this example consists
of only one entry) has negative signs. The second matrix $A_5$ has its first
hook negative, the second hook positive. The third matrix $A_6$ repeats the
behavior of the first one, the fourth one $A_7$ imitates the second one, and
so on.

Extensive experimentation shows that this {\em double alternating property}
holds for values of $l$ greater than zero. For coefficient matrices in the
traditional expansion range, the first matrix has its first hook positive,
the second one negative, the third one positive, etc. The second matrix has
the same alternating pattern of signs for the hooks but its first hook is
negative. The third matrix imitates the first, etc.

The following picture captures the phenomenon described above
for $n$ larger than one and when
the index $k$ is in the traditional range.
\[
\begin{array}{ccccc}
+ & + & + & \dots & + \\
+ & - & - & \dots & - \\
+ & - & + & \dots & + \\
+ & - & + \\
\vdots & \vdots & \vdots \\
+ & - & +
\end{array} \hskip 0.5 in \begin{array}{ccccc}
- & - & - & \dots & - \\
- & + & + & \dots & + \\
- & + & - & \dots & - \\
- & + & - \\
\vdots & \vdots & \vdots \\
- & + & -
\end{array} \hskip 0.5 in \mbox{etc.}
\]

\section{The relation with matrix valued orthogonal polynomials}

We close the paper remarking, once again, that our matrix valued spherical
functions are orthogonal with respect to a nice inner product and have
polynomial entries. Yet, they do not fit directly into the existing theory
of matrix valued orthogonal polynomials as given for instance in \cite{D} and
\cite{DVA}.

It is however possible to establish such a connection: define
the matrix valued function
$\Psi(j,t)$ by means of the relation

\[
\Psi(j,t) = \Phi(j,t)  \Phi^{-1}(0,t).
\]

It is now a direct consequence of the definitions that the family $\Psi(j,t)$
satisfies all the standard requirements in \cite{DVA} and not only satisfies
a three term recursion relation but also $\Psi(j,t)^t$ satisfies
a fixed differential equation with
matrix coefficients and only the "eigenvalue matrix"depends on $j$. In other
words the family $\Psi(j,t)$ meets all the conditions given at the beginning
of section $3$ and meets also the conditions of the standard theory in 
\cite {DVA} giving an example of a {\em classical family of matrix valued
orthogonal polynomials}. In particular, the coeffients in the differential
operator $D$ (obtained by conjugation from the one in \cite{GPT1} )
are matrix polynomials of degree going with the order of
differentiation.
For a nice introduction to this circle of ideas,
see the pioneering work in \cite{D}. 

In conclusion, we are really indebted to the editors for suggesting a number
of places where the exposition could be made a bit less terse. One of us, A.G.,
acknowledges a useful conversation with A. Duran that steered him in the
direction to section $5$ above.

\newcommand\bibit[5]{\bibitem
{#2}#3, {\em #4,\!\! } #5}

\end{document}